\documentclass[11pt,letterpaper]{article}

\usepackage[margin=1.1in]{geometry}
\usepackage{amsmath,amssymb,amsthm,mathtools}
\usepackage{enumitem}

\usepackage[colorlinks=true,
            linkcolor=blue,
            citecolor=blue,
            urlcolor=blue]{hyperref}

\newtheorem{theorem}{Theorem}[section]
\newtheorem{lemma}[theorem]{Lemma}
\newtheorem{proposition}[theorem]{Proposition}
\newtheorem{corollary}[theorem]{Corollary}
\theoremstyle{definition}
\newtheorem{definition}[theorem]{Definition}
\theoremstyle{remark}
\newtheorem{remark}[theorem]{Remark}

\newcommand{\R}{\mathbb R}
\newcommand{\Pp}{\mathbb P}
\newcommand{\Ee}{\mathbb E}
\newcommand{\cS}{\mathcal S}
\newcommand{\cG}{\mathcal G}
\newcommand{\cE}{\mathcal E}

\newcommand{\op}{\mathrm{op}}
\newcommand{\dens}{\operatorname{dens}}
\newcommand{\dist}{\operatorname{dist}}
\newcommand{\rank}{\operatorname{rank}}

\newcommand{\dpr}{d_{\mathrm{pr}}}
\newcommand{\1}{\mathbf 1}
\newcommand{\eps}{\varepsilon}

\title{Universality of smallest singular values of random row submatrices and phase retrieval problems}

\author{
  Shu Wen\\
  {\small University of Science and Technology of China}\\
  {\small\texttt{wenshu433@gmail.com}}
  \and
  Christian E. H\"aggblom\\
  {\small\texttt{chris.haggblo@gmail.com}}
  \and
  Xiufan Yang\\
  {\small Nanjing University of Posts and Telecommunications
}\\
  {\small\texttt{yxf3173616612@gmail.com}}
}

\date{June 28, 2026}

\begin{document}

\maketitle

\begin{abstract}
We establish universal logarithmic laws for extreme least singular values of random row submatrices and for the global stability of real phase retrieval near its injectivity threshold. Let $A_m\in\R^{N_m\times m}$ have i.i.d. entries distributed as a centered, variance-one real subgaussian random variable with a bounded Lebesgue density, and suppose $N_m/m\to\gamma>1$. For integers $0\le r_m\le N_m-m$, we study the minimum least singular value over the dependent family of row submatrices with $m+r_m$ rows,
\[
 M_{m,r_m}(A_m):=
 \min_{\substack{T\subset[N_m]\\|T|=m+r_m}}
 \sigma_m((A_m)_T).
\]
Writing $k_m=r_m+1$, we prove that, whenever $k_m\log m=o(m)$,
\[
 \frac{k_m}{m}\log M_{m,r_m}(A_m)
 \xrightarrow{\Pp}-h(\gamma),
 \qquad
 h(\gamma)=\gamma\log\gamma-(\gamma-1)\log(\gamma-1).
\]
This covers both fixed and moderately growing surplus. The limiting rate is universal within the stated class of entry distributions.

For real phase retrieval from $N_m=2m-1+r_m$ amplitude measurements, let $L_{A_m}$ denote the optimal lower Lipschitz constant of the map $x\mapsto|A_mx|$, with signals identified up to sign, and let $\kappa_{\mathrm{ph}}(A_m)=\|A_m\|_{\op}/L_{A_m}$ be its scale-invariant global condition number. In the same surplus regime, we establish
\[
 \frac{k_m}{m}\log L_{A_m}\xrightarrow{\Pp}-\log4,
 \qquad
 \frac{k_m}{m}\log\kappa_{\mathrm{ph}}(A_m)
 \xrightarrow{\Pp}\log4.
\]
Although these measurements determine every signal up to sign almost surely, the condition number grows exponentially in dimension for fixed surplus and exceeds every fixed power of $m$ throughout the stated regime, with probability tending to one. These conclusions quantify the distinction between uniqueness and uniform stability under measurement noise. All stated limits also hold almost surely on any common probability space, without requiring independence across dimensions.
\end{abstract}

\tableofcontents

\section{Introduction and main results}
Let $A\in\R^{N\times m}$, with $N\ge m$. For an integer $0\le r\le N-m$, define
\begin{equation}\label{eq:def-M}
 M_{m,r}(A):=
 \min_{\substack{T\subset[N]\\ |T|=m+r}}\sigma_m(A_T),
 \qquad M_m(A):=M_{m,0}(A),
\end{equation}
where $A_T$ is the row submatrix indexed by $T$. For any matrix $B$ with
$m$ columns, we use the convention
\[
 \sigma_m(B):=\inf_{\|x\|_2=1}\|Bx\|_2
\]
to denote  it's least singular value. In particular, $\sigma_m(B)=0$ when $B$ has rank less than $m$, including
when it has fewer than $m$ rows. The quantity $M_{m,r}(A)$ measures the
smallest least singular value among the $(m+r)\times m$ row submatrices;
the square case is $r=0$.

For background on least singular values of i.i.d. random matrices, see
\cite{RudelsonVershynin2009,TaoVu2010}. Even when the entries of $A$ are
independent, the family in \eqref{eq:def-M} contains exponentially many
strongly dependent matrices. A one-matrix estimate alone therefore does
not settle the problem. Related exponential conditioning questions occur
in smoothed analyses of random polytopes and Frank--Wolfe methods;
see \cite{RademacherShu}.

For $\gamma>1$, set
\begin{equation}\label{eq:def-h}
h(\gamma):=\gamma\log\gamma-(\gamma-1)\log(\gamma-1).
\end{equation}
If $N_m/m\to\gamma$, then
\[
\binom{N_m}{m}=\exp\{h(\gamma)m+o(m)\}.
\]
Shmalo \cite{ShmaloSquare} established the square-row law for real
standard Gaussian entries,
\[
 \frac1m\log M_m(A_m)\xrightarrow{\Pp}-h(\gamma),
\]
and the complex Gaussian exponent $-h(\gamma)/2$.
The corresponding real Gaussian results with surplus rows are developed
in \cite{ShmaloPrivate}. The lower-tail argument extends to
bounded-density entries by projection-density estimates. Obtaining the
matching upper bound requires control of exponentially small
small-ball probabilities over a strongly dependent family of random
hyperplanes. Our main theorems establish the real bounded-density
subgaussian law for both fixed and moderately growing surplus.

\subsection{Main theorems}

Throughout the paper, $\xi$ is a centered, variance-one real subgaussian
random variable with a bounded Lebesgue density $f_\xi$. Write
\[
 K:=\|f_\xi\|_\infty<\infty.
\]
All matrix entries are independent copies of $\xi$. No smoothness of
$f_\xi$ is assumed.

\begin{theorem}[Fixed-surplus law]\label{thm:main}
Let $r\geq0$ be a fixed integer, and let $A_m\in\R^{N_m\times m}$ have independent entries distributed as $\xi$. Assume
$N_m/m\to\gamma>1$ and $N_m\geq m+r$ for all large $m$. Then
\begin{equation}\label{eq:fixed}
\frac1m\log M_{m,r}(A_m)\xrightarrow{\Pp}-\frac{h(\gamma)}{r+1}
\end{equation}
\end{theorem}
\begin{remark}
Note that the normalization $\mathrm{Var}(\xi)=1$ is adopted only for notational convenience. In fact, the theorem remains valid for any fixed positive variance. The reasons are as follows: assume $\operatorname{Var}(\xi)=\sigma^2$, since $\xi$ is subgaussian and has a Lebesgue density, $0<\sigma<\infty$. Define
$
\widetilde{\xi}=\xi/\sigma,
\widetilde{A}_m=A_m/\sigma
$.
Then $\widetilde{\xi}$ is centered, variance one, subgaussian, and has density
$
f_{\widetilde{\xi}}(x)=\sigma f_{\xi}(\sigma x),
$
so that
$
\|f_{\widetilde{\xi}}\|_{\infty}\leq \sigma K.
$
Moreover,
$
M_{m,r}(A_m)=\sigma M_{m,r}(\widetilde{A}_m).
$
We have
$$
\frac{1}{m}\log M_{m,r}(A_m)
=
\frac{\log\sigma}{m}
+
\frac{1}{m}\log M_{m,r}(\widetilde{A}_m).
$$
Now since for every fixed $\varepsilon\in(0,h(\gamma))$, $|\log\sigma|/m\leq\varepsilon/2$ for sufficiently large $m$,
\begin{align*}
    \mathbb P\left\{\left|\frac{1}{m}\log M_{m,r}(A_m)+\frac{h(\gamma)}{r+1}\right|>\varepsilon\right\}&=\mathbb P\left\{\left|\frac{\log\sigma}{m}+\frac{1}{m}\log M_{m,r}(\tilde A_m)+\frac{h(\gamma)}{r+1}\right|>\varepsilon\right\}\\
    &\leq\mathbb P\left\{\left|\frac{1}{m}\log M_{m,r}(\tilde A_m)+\frac{h(\gamma)}{r+1}\right|>\frac{\varepsilon}{2}\right\}
\end{align*}
Thus, it is enough to establish the desired estimate for the normalized matrix $\widetilde{A}_m$, i.e., under the assumption $\mathrm{Var}(\xi)=1$. For the remainder of the proof we impose this normalization.
\end{remark}
If $r$ is not fixed, we have the following moderate-surplus form.
\begin{theorem}[Moderate-surplus law]\label{thm:moderate}
Let $N_m$ and $r_m$ be integer-valued sequences, and let $A_m\in\R^{N_m\times m}$ have i.i.d. entries distributed as $\xi$. Suppose
\[
 \frac{N_m}{m}\longrightarrow\gamma>1,
 \qquad
 k_m=r_m+1\in\mathbb N,
 \qquad
 k_m\log m=o(m),\qquad0\le r_m\le N_m-m
\]
for all large $m$. Then
\begin{equation}\label{eq:th1.1}
 \frac{k_m}{m}\log M_{m,r_m}(A_m)
 \xrightarrow{\Pp}-h(\gamma).
\end{equation}
\end{theorem}

The deviation estimates proved below are summable. Consequently, if
the matrices $A_m$ are defined on a common probability space, the
limits in Theorems~\ref{thm:main} and~\ref{thm:moderate} hold almost
surely, without any independence assumption across dimensions;
see Corollary~\ref{cor:almost-sure}. In the phase-retrieval setting,
the same argument gives almost-sure logarithmic laws for the
obstruction $\omega(A_m)$ and the optimal lower Lipschitz constant
$L_{A_m}$, together with the corresponding condition-number rate
in Corollary~\ref{cor:phase-condition-rate}.

The scale in Theorem~\ref{thm:moderate} comes from balancing the number
of candidate row sets against the probability of satisfying the smallness
constraints. Write $k=k_m=r_m+1$. There are
$\exp\{h(\gamma)m+o(m)\}$ candidate row sets. An $(m-1)$-row base set
determines a unit normal $u$, and adjoining $k$ rows whose scalar products
with $u$ have absolute value at most $t/\sqrt{k}$ produces an
$(m+k-1)\times m$ matrix with least singular value at most $t$.
For a delocalized normal, each independent constraint has probability
comparable to $t/\sqrt{k}$. Their combined cost is therefore
$t^k\exp\{O(k\log(k+1))\}$. Since $k\log m=o(m)$, balancing this cost
against the number of candidates suggests
$k\log t\sim-h(\gamma)m$. The upper-bound proof makes this heuristic
rigorous despite dependence, by reserving a small proportion of the rows
and imposing the constraints sequentially while controlling the
surviving pairs of candidate normals.

\begin{remark}[Scope of the assumptions]\label{rem:scope}
The assumptions on $\xi$ are sufficient and are not asserted to be
optimal. Absolute continuity of the joint entry distribution implies
that every prescribed row or column submatrix has its maximal possible
rank almost surely: rank deficiency forces a nonzero minor polynomial
to vanish, and its zero set has Lebesgue measure zero. Taking a finite
union gives the simultaneous rank properties needed at each dimension.
Thus the logarithms and normal lines used below are well defined without
any smoothness assumption.

The bounded-density hypothesis supplies the projection-density bounds
in Theorem~\ref{thm:projection-density}, including the one- and
two-dimensional anti-concentration estimates, and the density bounds
for the blocks used in Lemma~\ref{lem:local-lower}. The density-level
Berry--Esseen estimate itself requires bounded densities and the stated
variance and third-moment bounds, not differentiability of the entry
density. Subgaussianity provides those moment bounds, the operator-norm
tails, and the normal-vector delocalization estimate of
Theorem~\ref{thm:NV}. Atomic laws fall outside these assumptions and
require separate treatment; some such laws also have exact
phase-retrieval obstructions.
\end{remark}

\subsection{Application to phase retrieval}
\label{seubsec:intro-phase-retrieval}

 The row-submatrix problem also arises naturally in real phase retrieval, where one seeks to recover a signal from the magnitudes of its linear measurements. For a measurement matrix $A\in\R^{N\times m}$, the amplitude map is
 \[
    \Phi_A:\R^m/\{\pm1\}\longrightarrow\R^N
    \qquad \Phi_A(x)=|Ax|
 \]
with the absolute value taken coordinatewise. The quotient accounts for the unavoidable ambiguity between $x$ and $-x$. Besides the uniqueness of recovery, a basic question is how sensitively the recovered signal can depend on perturbations of the measured amplitudes.

The connection with row submatrices comes from the loss of sign information. When two signals $x$ and $y$ have similar measurement magnitudes, their signed measurements can be close either with the same sign or with opposite signs. The measurement rows can therefore be partitioned into a set $J$, on which $A_J(x-y)$ is small, and its complement, on which $A_{J^c}(x+y)$ is small. If the signals are separated up to sign, both $x-y$ and $x+y$ are bounded away from zero. Such a pair can consequently occur only when both complementary row blocks have small least singular values. Conversely, directions that are poorly detected by complementary blocks can be combined to produce signals with similar amplitude measurements. This connects phase-retrieval stability to the extreme singular values of row submatrices. For random measurement matrices, our estimates quantify the resulting instability near the critical measurement count.

Balan, Casazza, and Edidin~\cite{BalanCasazzaEdidin} characterized injectivity in the real case by the complement property: for every $J\subset[N]$, either the rows of $A_J$ or those of $A_{J^c}$ span $R^m$. Consequently, $2m-1$ is the sharp generic measurement threshold for recovery of every real signal up to sign. Below this count no measurement family has the complement property, whereas at or above it every full-spark matrix does. Here full-spark means that every $m$ rows are linearly independent. This holds almost surely under our absolute-continueity assumption.

Injectivity at each dimension leaves open the quantitative behavior of the inverse as the dimension grows. With $d_\pm(x,y)=\min\{\|x-y\|_2,\|x+y\|_2\}$, the optimal lower Lipschitz constant and the associated condition number are
\[
    L_A:=\inf_{x\ne\pm y} \frac{\|\Phi_A(x)-\Phi_A(y)\|_2}{d_\pm(x,y)}
    \qquad
    \kappa_{\mathrm{ph}}(A):=\frac{\|A\|_{\op}}{L_A}
\]
where the latter is defined when $L_A>0$. The optimal upper Lipschitz constant is $\|A\|_{\op}$. Thus a small $L_A$ means that some separated signals have nearly indistinguishable amplitude measurements, while $\kappa_{\mathrm{ph}}(A)$ measures global conditioning independently of a common rescaling of the measurements.

Quantitative versions of the complement property connect this stability problem to singular values of complementary row blocks. Balan and Wang~\cite{BalanWang} identified the optimal lower Lipschitz constant through an exact partition formula. Bandeira, Cahill, Mixon, and Nelson\cite{BandeiraCahillMixonNelson} introduced the strong complement property and established quantitative stability estimates for Gaussian measurements. In particular, for each fixed measurement-to-dimension ratio strictly greater than two, their estimates give a condition number bounded independently of dimension with probability tending to one. This raises the question of how conditioning deteriorates when the measurement count approaches the injectivity threshold.

The connection with our row minima is particularly direct when $N=2m-1+r$ and $A$ is full spark. The Balan--Wang obstruction then satisfies
\[
    \omega(A):=
    \min_{\substack{J\subset[N]\\\rank(A_J)<m}}
        \sigma_m(A_{J^c})
    =M_{m,r}(A), \quad L_A\le\omega(A)
\]
At the critical count $r=0$, every partition has a rank-deficient side, and the paritition formula yields $L_A=\omega(A)=M_{m,0}(A)$. With additional rows, both sides can have full rank, so identifying the rate of $L_A$ requires further estimate; see Proposition~\ref{prop:phase-partition} and Lemma~\ref{lem:omega-reduction}.

For Gaussian entries and fixed surplus $r$, the results in \cite{ShmaloSquare,ShmaloPrivate} identify the exponential decay of the row minimum through $m^{-1}\log M_{m,r}(A_m)\xrightarrow{\Pp}-h(\gamma)/(r+1)$ when $N_m/m\to\gamma>1$, yielding the rate $-\log4$ at the critical phase-retrieval count $N_m=2m-1$. At the critical count, H\"aggblom~\cite{HaggblomGaussian} obtained finer localization of $L_A=\omega(A)$ around the scale $\bigl(\sqrt m\binom{2m-1}{m}\bigr)^{-1}$, up to arbitrary diverging subexponential factors. The present paper establishes logarithmic universality for bounded-density subgaussian entries and treats both fixed and moderately growing surplus.

More precisely, let $r_m\ge0$ be intergers, let $A_m$ have entries satisfying the standing assumptions, and suppose
\[
    N_m=2m-1+r_m,\qquad k_m=r_m+1, \qquad k_m\log m=o(m)
\]
Then $N_m/m\to2$ and $h(2)=\log4$, so Theorem~\ref{thm:moderate} determines the logarithmic rate of $\omega(A_m)$. Corollaries~\ref{cor:phase-lower-Lipschitz} and \ref{cor:phase-condition-rate} establish
\[
    \frac{k_m}{m}\log L_{A_m}\xrightarrow{\Pp}-\log4
    \qquad
    \frac{k_m}{m}\log\kappa_{\mathrm{ph}}(A_m)
    \xrightarrow{\Pp}\log4
\]
The lower estimate for $L_{A_m}$ combines the partition formula with independent small-ball estimates for complementary blocks and a union bound, following the strategy of the Gaussian strong-complement argument in~\cite[Theorem~20]{BandeiraCahillMixonNelson}. The matching upper estimate follows from $L_{A_m}\le\omega(A_m)$ and our row-minimum theorem. These limits also hold almost surely on any common probability space, without independence across dimensions, by Corollaries~\ref{cor:almost-sure} and \ref{cor:phase-condition-rate}.

These conclusions distinguish exact recovery from robust recovery. Under our assumptions, the measurements determine every signal up to sign almost surely at each dimension. Nevertheless, some distinct signals produce measurement vectors that are extremely close together. To interpret this independently of the overall measurement scale, normalize $A_m$ by $\|A_m\|_{\op}$. The resulting amplitude map has optimal upper Lipschitz constant one and optimal lower Lipschitz constant $1/\kappa_{\mathrm{ph}}(A_m)$. Thus there are signal pairs separated by unit distance, modulo sign, whose normalized measurement vectors are separated by only $1/\kappa_{\mathrm{ph}}(A_m)$. Measurement perturbations of this order can make such signals indistinguishable. This difficulty is intrinsic to the measurements. Any reconstruction method that recovers all signals must contend with the corresponding worst-case sensitivity.

For fixed surplus $r$, this sensitivity grows exponentially with dimension, at logarithmic rate $\log4/(r+1)$. Additional measurements therefore reduce the rate of deterioration, but a fixed number of rows beyond the injectivity threshold does not prevent exponential ill-conditioning. The conclusion also covers growing surplus. Whenever $(r_m+1)\log m=o(m)$, the condition number exceeds every fixed power of $m$ with probability tending to one. In particular, the additional measurements in this regime are insufficient to obtain even a polynomial bound on the global condition number, despite guaranteeing uniqueness of recovery.

The result therefore quantifies the gap between uniqueness and uniform robustness near the critical measurement count. Its interpretation is a worst-case one, concerning recovery uniformly over all signals. Particular signals or restricted signal classes may have substantially better stability. The common logarithmic rate of $L_{A_m}$ and $\omega(A_m)$ identifies the scale of this global instability, without asserting equality of these quantities for redundant matrices or a precise multiplicative asymptotic.

\subsection{Comparison with the Gaussian argument}

The Gaussian proof uses rotational invariance to describe hyperplane
normals by Haar-distributed directions. For general i.i.d. entries,
we instead use the following ingredients.
\begin{enumerate}[label=\textup{(\roman*)},leftmargin=2.5em]
\item A density-level local central limit estimate gives a uniform
positive density near zero for $\langle X,u\rangle$ when $u$ is
sufficiently delocalized.
\item The Nguyen--Vu normal-vector estimate, followed by averaging
and Markov's inequality, shows that a positive fraction of the
candidate hyperplanes have delocalized normals.
\item Except for an exponentially small fraction of pairs, the normals
have a fixed positive projective separation. This follows from an
operator-norm bound and fresh-row small-ball estimates.
\item The reserved rows are split into $k_m$ independent groups.
A one-row second-moment argument provides a fixed positive success
probability; selection within each group amplifies this probability.
The surviving candidates and exceptional pairs are controlled through
all $k_m$ restrictions.
\end{enumerate}

\section{Preliminaries}\label{sec:preliminaries}
\subsection{Notation}
We use the usual small $o$ and big $O$ notation. For a real random variable $Y$, let $\dens(Y)$ denote its
Lebesgue density and let $\|\dens(Y)\|_\infty$ denote the essential supremum of that density. We
use the Orlicz convention
\[
\|Y\|_{\psi_2}:=\inf\bigl\{s>0:\Ee\exp(Y^2/s^2)\le2\bigr\},
\]
and say that $Y$ is subgaussian when $\|Y\|_{\psi_2}<\infty$.

Throughout the proof, constants denoted by $c,C,C_1\dots$ may change from line to line.
Constants in distributional lemmas depend only on the law of $\xi$; constants in tail estimates may
also depend on fixed parameters such as $\gamma$ and $\varepsilon$, but never on $m$.

For unit vectors $u,v\in\R^m$, define their projective distance by
\begin{equation}\label{eq:projective-distance}
\dpr(u,v):=\min\{\|u-v\|_2,\|u+v\|_2\}.
\end{equation}
\subsection{Upper bounds for projection densities}

We first record the projection-density theorem of Rudelson and Vershynin
\cite{RudelsonVershyninDensity}. The distinction between its general $d$-dimensional form and its
sharp one-dimensional form is important for the two-dimensional estimate below.

\begin{theorem}[Densities of projections]\label{thm:projection-density}
Let $X=(X_1,\dots,X_n)$ have independent real-valued coordinates whose densities are bounded by
$K$ almost everywhere. If $P$ is the orthogonal projection in $\R^n$ onto a $d$-dimensional
subspace, then the density of $PX$, with respect to $d$-dimensional Lebesgue measure on that
subspace, is bounded by $(C_{\mathrm{RV}}K)^d$ almost everywhere, where
$C_{\mathrm{RV}}>0$ is an absolute constant. In dimension $d=1$, the sharper bound
$\sqrt2K$ holds.
\end{theorem}

\begin{lemma}[One- and two-dimensional marginal bounds]\label{lem:marginal-density}
Let $\xi$ be as in Theorem~\ref{thm:main}, and let $X=(\xi_1,\dots,\xi_m)$ have independent
coordinates distributed as $\xi$. Then, for every deterministic unit vector $u\in\R^m$,
\begin{equation}\label{eq:one-dimensional-density}
\|\dens(\langle X,u\rangle)\|_\infty\le\sqrt2K.
\end{equation}
Furthermore, for unit vectors $u,v\in\R^m$ with $\rho=\langle u,v\rangle$ and $|\rho|<1$,
\begin{equation}\label{eq:two-dimensional-density}
\|\dens(\langle X,u\rangle,\langle X,v\rangle)\|_\infty
\le \frac{C_{\mathrm{RV}}^2K^2}{\sqrt{1-\rho^2}}.
\end{equation}
\end{lemma}

\begin{proof}
Let $E=\operatorname{span}(u)$ and let $P_E$ be the orthogonal projection onto $E$. Since
$P_EX=\langle X,u\rangle u$, the isometry $t\mapsto tu$ identifies the density of $P_EX$ on $E$
with the density of $\langle X,u\rangle$. The sharp one-dimensional part of
Theorem~\ref{thm:projection-density} therefore gives \eqref{eq:one-dimensional-density}.

For the two-dimensional estimate, let $E=\operatorname{span}(u,v)$ and choose the orthonormal
basis
\[
e_1=u,
\qquad
e_2=\frac{v-\rho u}{\sqrt{1-\rho^2}}.
\]
Set $Y=(\langle X,e_1\rangle,\langle X,e_2\rangle)$. The general $d=2$ part of
Theorem~\ref{thm:projection-density} gives
\[
\|\dens(Y)\|_\infty\le C_{\mathrm{RV}}^2K^2.
\]
Moreover,
\[
(\langle X,u\rangle,\langle X,v\rangle)=BY,
\qquad
B=\begin{pmatrix}1&0\\ \rho&\sqrt{1-\rho^2}\end{pmatrix},
\]
and $|\det B|=\sqrt{1-\rho^2}$. The change-of-variables formula for densities yields
\[
\|\dens(BY)\|_\infty
\le \frac{\|\dens(Y)\|_\infty}{|\det B|},
\]
which is exactly \eqref{eq:two-dimensional-density}.
\end{proof}

We next extract the one-dimensional form of the density-level Berry--Esseen theorem needed below.

\begin{theorem}[Density-level Berry--Esseen theorem]\label{thm:density-BE}
Let $Y_1,\dots,Y_r$ be independent centered real random variables such that
\[
\sum_{\ell=1}^r\operatorname{Var}(Y_\ell)=r.
\]
Assume that, for fixed positive constants $b,M,B$, every $\ell$ satisfies
\[
\operatorname{Var}(Y_\ell)\le b,
\qquad
\|\dens(Y_\ell)\|_\infty\le M,
\qquad
\Ee|Y_\ell|^3\le B.
\]
If $f_r$ is the density of $r^{-1/2}\sum_{\ell=1}^rY_\ell$ and $\varphi$ is the standard Gaussian
density, then
\begin{equation}\label{eq:density-BE}
\sup_{x\in\R}|f_r(x)-\varphi(x)|\le\frac{C(b,M,B)}{\sqrt r}.
\end{equation}
\end{theorem}

\begin{proof}
Apply the one-dimensional case of
\cite[Theorem~1, inequality~(1.1)]{BobkovGoetze}, with
$n=r$, $X_\ell=Y_\ell$, and
$Z_n=r^{-1/2}\sum_{\ell=1}^rY_\ell$. In that theorem the parameters are
\[
 \sigma^2=\max_\ell\operatorname{Var}(Y_\ell),\qquad
 M_*=\max_\ell\|\dens(Y_\ell)\|_\infty,\qquad
 \beta_3=\frac1r\sum_{\ell=1}^r\Ee|Y_\ell|^3.
\]
The variance normalization here is precisely the normalization there.
Consequently,
\[
 \|f_r-\varphi\|_\infty
 \le \frac{C\sigma M_*^2\beta_3}{\sqrt r}
 \le \frac{C\sqrt b\,M^2B}{\sqrt r},
\]
which proves \eqref{eq:density-BE}. No smoothness of the densities is
required.
\end{proof}

The advantage of this density-level estimate is that its integrated error
scales with the length of the interval. If
$\|f-\varphi\|_\infty\le\delta$, then for every $t>0$,
\[
 \left|\int_{-t}^t f(x)\,dx-\int_{-t}^t\varphi(x)\,dx\right|
 \le 2t\delta.
\]
Thus a sufficiently small uniform density error gives a useful relative
estimate even when $t$ is exponentially small. An approximation of
distribution functions with uniform error $\delta$ gives only an
absolute interval error at most $2\delta$, which may exceed the
probability being estimated.

Combining Lemma~\ref{lem:marginal-density} with Theorem~\ref{thm:density-BE} gives the
local lower bound used in the upper-tail argument. The point is that, when $u$ is delocalized, no
single summand $u_j\xi_j$ contributes appreciably to $\langle u,X\rangle$. The resulting sum lies
in a local central-limit regime, so its density is uniformly positive near the origin.
\subsection{A local lower bound for projection densities}

\begin{lemma}[Uniform local lower density]\label{lem:local-lower}
There are constants $\eta,c_0>0$, depending only on the law of $\xi$, such that for every $m$ and every unit vector $u\in\R^m$ with $\|u\|_\infty\le\eta$, a version $f_u$ of the density of $\langle X,u\rangle$ satisfies
\begin{equation}\label{eq:local-density-lower}
 f_u(x)\ge c_0,\quad\forall|x|\le1.
 \end{equation}
\end{lemma}
\begin{proof}
Fix an integer $q\ge1$, to be chosen sufficiently large in terms of the distribution of $\xi$. We
first assume that
\[
\max_{1\le j\le m}u_j^2\le\frac1{4q}.
\]
Set $w_j:=u_j^2$. Since $\|u\|_2=1$, we have $\sum_{j=1}^mw_j=1$.

We partition the nonzero weights $w_j$ into consecutive blocks by the following greedy procedure.
Starting with an empty block, add nonzero weights until their sum first becomes at least $1/q$, then
close the block and repeat with the remaining weights. Since each weight is at most $1/(4q)$, every
completed block has total weight between $1/q$ and $5/(4q)$. At the end there may remain a
terminal collection of total mass less than $1/q$; merge it into the last completed block.
Because the total mass is one, at least one block is completed. We obtain nonempty disjoint
blocks $B_1,\dots,B_r$ whose union contains all $j$ with $u_j\ne0$, and, with
$w_\ell:=\sum_{j\in B_\ell}u_j^2$,
\[
\frac1q\le w_\ell\le\frac9{4q},
\qquad
\sum_{\ell=1}^rw_\ell=1.
\]
The upper bound follows because the last completed block has weight at most $5/(4q)$ before the
terminal remainder, while the remainder has weight less than $1/q$. Consequently,
\[
\frac{4q}{9}\le r\le q.
\]

For $\ell=1,\dots,r$, define
\[
Z_\ell:=\sqrt r\sum_{j\in B_\ell}u_j\xi_j.
\]
The variables $Z_1,\dots,Z_r$ are independent and centered. Since
$\operatorname{Var}(\xi)=1$,
\[
\operatorname{Var}(Z_\ell)=r\sum_{j\in B_\ell}u_j^2=rw_\ell,
\]
and therefore
\[
\frac49\le\operatorname{Var}(Z_\ell)\le\frac94,
\qquad
\sum_{\ell=1}^r\operatorname{Var}(Z_\ell)=r.
\]

We next establish a uniform density bound. If $B_\ell$ contains at least two indices, put
$a_j:=u_j/\sqrt{w_\ell}$ for $j\in B_\ell$. Then $\sum_{j\in B_\ell}a_j^2=1$, so
Lemma~\ref{lem:marginal-density} gives
\[
\left\|\dens\left(\sum_{j\in B_\ell}a_j\xi_j\right)\right\|_\infty\le\sqrt2K.
\]
Since
\[
Z_\ell=\sqrt{rw_\ell}\sum_{j\in B_\ell}a_j\xi_j,
\]
the scaling rule for densities yields
\[
\|\dens(Z_\ell)\|_\infty
\le\frac{\sqrt2K}{\sqrt{rw_\ell}}
\le\frac{3\sqrt2}{2}K.
\]
If $B_\ell=\{j\}$ is a singleton, then $Z_\ell=\sqrt r\,u_j\xi_j$, and
\[
\|\dens(Z_\ell)\|_\infty
\le\frac{K}{\sqrt r\,|u_j|}
=\frac{K}{\sqrt{rw_\ell}}
\le\frac32K.
\]
Thus $\|\dens(Z_\ell)\|_\infty\le C_1K$ in all cases.

The third moments are also uniformly bounded. Since $\xi$ is subgaussian and the coordinates are
independent and centered, the standard subgaussian inequality for linear combinations
\cite[Proposition~2.6.1]{VershyninBook} gives
\[
\left\|\sum_{j\in B_\ell}u_j\xi_j\right\|_{\psi_2}
\le C\|\xi\|_{\psi_2}\left(\sum_{j\in B_\ell}u_j^2\right)^{1/2}.
\]
Consequently,
\[
\|Z_\ell\|_{\psi_2}
\le C\|\xi\|_{\psi_2}\sqrt{rw_\ell}
\le\frac{3C}{2}\|\xi\|_{\psi_2},
\]
so $\Ee|Z_\ell|^3\le C_2$, where $C_2$ depends only on the law of $\xi$.

We may therefore apply Theorem~\ref{thm:density-BE} with
$b=9/4$, $M=C_1K$, and $B=C_2$. Since
\[
\frac1{\sqrt r}\sum_{\ell=1}^rZ_\ell
=\sum_{j=1}^mu_j\xi_j
=\langle X,u\rangle,
\]
if $f_u$ denotes the density of $\langle X,u\rangle$, then
\[
\sup_{x\in\R}|f_u(x)-\varphi(x)|
\le\frac{C_3}{\sqrt r}
\le\frac{3C_3}{2\sqrt q}.
\]
Choose $q$ sufficiently large that $3C_3/(2\sqrt q)\le\varphi(1)/2$. Then, for $|x|\le1$,
\[
f_u(x)\ge\varphi(x)-\frac{\varphi(1)}2\ge\frac{\varphi(1)}2.
\]
Set
\[
\eta:=\frac1{2\sqrt q},
\qquad
c_0:=\frac{\varphi(1)}2.
\]
If $\|u\|_\infty\le\eta$, then $\max_ju_j^2\le1/(4q)$, so all preceding estimates apply and we get \eqref{eq:local-density-lower}.
\end{proof}

\section{The lower estimate}

The lower estimate follows from a small-ball bound for one rectangular
matrix and a union bound over row sets. We first record the required
count and then state the one-matrix estimate separately.

\begin{lemma}[Binomial estimates]\label{lem:binomial}
Let $N_m/m\to\gamma>1$, let $k_m\ge1$ be integers, and assume
$k_m\log m=o(m)$. Then
\[
 \log\binom{N_m}{m+k_m-1}=h(\gamma)m+o(m).
\]
\end{lemma}
\begin{proof}
Set $a_m=(m+k_m-1)/m$ and $b_m=N_m/m$. Then $a_m\to1$ and
$b_m\to\gamma$. Stirling's formula gives
\[
 \frac1m\log\binom{N_m}{m+k_m-1}
 =b_m\log b_m-a_m\log a_m
  -(b_m-a_m)\log(b_m-a_m)+o(1)
 \longrightarrow h(\gamma).
\]
\end{proof}

\begin{lemma}[Rectangular small-ball estimate]\label{lem:rectangular-small-ball}
Let $B$ be an $(m+d-1)\times m$ matrix whose entries are independent
and have Lebesgue densities bounded by $K$, where $m,d\ge1$ are integers.
There is an absolute constant $C_{\mathrm{RV}}>0$ such that, for every
$t>0$,
\begin{equation}\label{eq:rectangular-small-ball}
 \Pp\{\sigma_m(B)\le t\}
 \le m(C_{\mathrm{RV}}K)^d v_d(\sqrt m\,t)^d
 \le m(D\sqrt m\,t)^d,
\end{equation}
where $v_d$ is the volume of the Euclidean unit ball in $\R^d$ and
$D:=\max\{1,2C_{\mathrm{RV}}K\}$ depends only on $K$.
\end{lemma}
\begin{proof}
Write $C_1,\ldots,C_m$ for the columns of $B$, and set
\[
 H_j:=\operatorname{span}\{C_\ell:\ell\ne j\},\qquad
 \delta_j:=\dist(C_j,H_j).
\]
Absolute continuity implies full column rank almost surely, both for
$B$ and for each matrix obtained by removing one column. The negative
second-moment identity is
\begin{equation}\label{eq:negative-second-moment}
 \sum_{i=1}^m\sigma_i(B)^{-2}
 =\operatorname{tr}\bigl((B^TB)^{-1}\bigr)
 =\sum_{j=1}^m\delta_j^{-2}.
\end{equation}
Indeed, the $j$th diagonal entry of $(B^TB)^{-1}$ is $\delta_j^{-2}$.
If $\sigma_m(B)\le t$, then some $\delta_j\le\sqrt m\,t$.

Condition on all columns except $C_j$. The subspace $H_j^\perp$
has dimension $d$ almost surely, and $C_j$ is independent of it.
Theorem~\ref{thm:projection-density} therefore gives
\[
 \Pp\{\delta_j\le\sqrt m\,t\mid (C_\ell)_{\ell\ne j}\}
 \le (C_{\mathrm{RV}}K)^d v_d(\sqrt m\,t)^d.
\]
Take expectations and sum over $j$. Finally, $v_d\le2^d$, since the
unit ball is contained in $[-1,1]^d$. This proves
\eqref{eq:rectangular-small-ball}.
\end{proof}

\begin{proposition}[Lower estimate]\label{prop:lower-tail}
Under the assumptions of Theorem~\ref{thm:moderate}, for every
$\varepsilon>0$,
\begin{equation}\label{eq:lower-tail-event}
 \Pp\left\{M_{m,r_m}(A_m)\le
 \exp\left[-(h(\gamma)+\varepsilon)\frac{m}{k_m}\right]\right\}
 \le \exp\{-\varepsilon m+o(m)\}\longrightarrow0.
\end{equation}
\end{proposition}
\begin{proof}
There are $\binom{N_m}{m+k_m-1}$ row submatrices of the required size.
By Lemma~\ref{lem:rectangular-small-ball} with $d=k_m$, for every $t>0$,
\[
 \Pp\{M_{m,r_m}(A_m)\le t\}
 \le \binom{N_m}{m+k_m-1}m(D\sqrt m\,t)^{k_m}.
\]
For $t=\exp[-(h(\gamma)+\varepsilon)m/k_m]$, the logarithm of the
right-hand side is
\[
 h(\gamma)m+o(m)+\log m+k_m\log D
 +\frac{k_m}{2}\log m-(h(\gamma)+\varepsilon)m
 =-\varepsilon m+o(m),
\]
using $k_m\log m=o(m)$. This proves the claim.
\end{proof}

\section{The upper estimate}
The upper estimate is substantially more delicate than the lower estimate.
Indeed, the lower estimate only requires ruling out the existence of an
exceptionally ill-conditioned row submatrix, and therefore follows from a
one-matrix small-ball bound combined with a union bound. In the opposite
direction, the upper estimate must prove that at least one such submatrix
actually exists. Although there are exponentially many candidate row
submatrices, they are strongly dependent because different candidates
typically share a large number of rows.

A further difficulty is that the Gaussian argument
\cite{ShmaloPrivate} uses rotational invariance of the normal vectors.
For general i.i.d. entries, those normals need not be Haar distributed.
The necessary small-ball estimates must instead hold uniformly, or
on average, over an exponentially large family of random hyperplanes,
at the exponentially small scale

\[
t=\exp\left(-\Theta\left(\frac{m}{k}\right)\right).
\]
Moreover, pairs of candidate normals may be nearly parallel, in which case
their small-ball events can be highly correlated. The proof must therefore control not only individual random normals but
also the dependence between different candidate hyperplanes.
\subsection{Proof map for the upper estimate.}
For convenience, write
\[
k=k_m=r_m+1,
\qquad
\lambda=h(\gamma)-\varepsilon,
\qquad
t=\exp\left(-\lambda\frac{m}{k}\right).
\]
Our goal is to prove that, with probability tending to one, there is a set
of $m+k-1=m+r_m$ rows whose least singular value is at most $t$.

We begin by choosing $kb$ rows and call them the marked rows, where
$
b=\left\lfloor\frac{\delta m}{k}\right\rfloor,
$
and divide them into $k$ disjoint groups of $b$ rows each. The remaining
\[
n_m:=N_m-kb
\]
rows will be called the base rows. Since $kb=\delta m+O(k)$ and $k=o(m)$,
\[
\frac{n_m}{m}\longrightarrow\beta:=\gamma-\delta.
\]
For every set $S$ of $m-1$ base rows, let $u_S$ be a unit vector orthogonal
to all rows indexed by $S$. The number of such sets is
\[
L=\binom{n_m}{m-1}
  =\exp\{h(\beta)m+o(m)\}.
\]
We choose $\delta>0$ sufficiently small that
$
h(\beta)>\lambda.
$
The reason for this choice will appear below: the initial number of sets
$S$ grows like $\exp\{h(\beta)m\}$, whereas the $k$ conditions imposed later
produce a main factor $\exp(-\lambda m)$. The strict inequality
$h(\beta)>\lambda$ ensures that some sets $S$ remain after all $k$
conditions have been imposed.

Section~4 establishes the estimates needed to compare the normal vectors
$u_S$. First, it shows that at least half of them have sufficiently small
coordinates. For such a vector $u_S$ and an independent row $X$, the
probability of
\[
|\langle X,u_S\rangle|\le s
\]
is bounded above and below by constant multiples of $s$. Geometrically,
this inequality means that $X$ is close to the hyperplane spanned by the
rows indexed by $S$. Second, Section~4 shows that only a very small
proportion of pairs $(S,T)$ have normal vectors $u_S$ and $u_T$ that are
close to each other, up to a change of sign. For all other pairs, the
probability that
\[
|\langle X,u_S\rangle|\le s
\quad\text{and}\quad
|\langle X,u_T\rangle|\le s
\]
hold simultaneously is bounded by a constant multiple of $s^2$. These
one-vector and two-vector estimates allow us to control both the number of
sets that remain and the dependence between different sets.

We then use the $k$ reserved groups one at a time. Set
$
s=t/\sqrt{k}.
$
In each round $j$, we examine the $b$ rows in the $j$-th group and choose one row
$X_j$. After choosing $X_j$, we keep only those sets $S$ for which
\[
|\langle X_j,u_S\rangle|\le s.
\]
Lemma~\ref{lem:one-row-stage} shows that each independent row has a fixed positive probability
of keeping sufficiently many sets $S$ while preserving the pairwise
estimates needed for the next round. Since each group contains
$b\asymp m/k$ independent rows, the probability that no row in a given
group works is at most
\[
\exp\left(-c\frac{m}{k}\right).
\]
Repeating this argument for all $k$ groups gives rows
$X_1,\ldots,X_k$ and at least one set $S$ satisfying
\[
|\langle X_j,u_S\rangle|\le s,
\qquad 1\le j\le k.
\]
The fact that at least one set remains follows from
$h(\beta)>\lambda$ and $k\log k=o(m)$.

Finally, form a matrix $B$ from the $m-1$ base rows indexed by $S$ and the
selected rows $X_1,\ldots,X_k$. The matrix $B$ has
\[
m-1+k=m+r_m
\]
rows. Since the base rows indexed by $S$ are orthogonal to $u_S$,
\[
\|Bu_S\|_2^2
=
\sum_{j=1}^k|\langle X_j,u_S\rangle|^2
\le ks^2=t^2.
\]
Because $\|u_S\|_2=1$, the variational characterization of the least
singular value gives
\[
\sigma_m(B)\le\|Bu_S\|_2\le t.
\]
Therefore
$
M_{m,r_m}(A_m)\le t.
$
The probability that the construction fails is bounded by
\[
Ce^{-cm}
+
k\exp\left(-c\frac{m}{k}\right),
\]
which tends to zero under the assumption $k\log m=o(m)$.
\subsection{Geometry of random hyperplane normals}\label{sec:geometry}

Because the entry law is absolutely continuous, with probability one every row family of size at
most $m$ that appears below has the expected rank. We work on this full-probability event and
choose each normal measurably, for example by requiring its first nonzero coordinate to be positive.

We use the largest-coordinate estimate
\cite[Theorem~1.4, equation~(6)]{NguyenVu}. The variables denoted by
$n$ and $m$ in that estimate correspond to $m$ and $s$, respectively,
in the formulation below.

\begin{theorem}[Nguyen--Vu]\label{thm:NV}
Let $B$ be an $(m-1)\times m$ matrix with i.i.d. entries distributed as
$\xi$, and let $u$ be a measurably chosen unit normal to its rows.
There are constants $C,C_1>0$, depending only on the entry law, such
that, for all sufficiently large $m$ and every $s\ge C_1\log m$,
\begin{equation}\label{eq:NV}
 \Pp\left\{\|u\|_\infty\ge\sqrt{\frac{s}{m}}\right\}
 \le Cm^2e^{-s/C}.
\end{equation}
\end{theorem}

Theorem~\ref{thm:NV} provides a delocalization estimate for the unit normal of a single random
hyperplane. The following lemma shows that, with high probability, at least half of the hyperplanes
generated by the random rows have delocalized normals.

\begin{lemma}[Many hyperplanes have delocalized normals]\label{lem:many-delocalized}

Let $n=n_m\ge m$, let $R_1,\dots,R_n$ be independent rows with i.i.d. coordinates
distributed as $\xi$, and for every $S\subset[n]$ with $|S|=m-1$ let $u_S$ be a unit normal to the
rows indexed by $S$. Put
\[
\cS:=\{S\subset[n]:|S|=m-1\},
\qquad
L:=|\cS|=\binom{n}{m-1},
\]
and, with $\eta$ from Lemma~\ref{lem:local-lower}, define the set of good hyperplanes
\[
\cG:=\{S\in\cS:\|u_S\|_\infty\le\eta\}.
\]
Then there are constants $C,c>0$ such that
\begin{equation}\label{eq:many-good}
\Pp\{|\cG|<L/2\}\le Cm^2e^{-cm}
\end{equation}
for all sufficiently large $m$.
\end{lemma}

\begin{proof}
Fix $S\in\cS$. The matrix whose rows are indexed by $S$ is an $(m-1)\times m$ random matrix
with independent entries distributed as $\xi$, so Theorem~\ref{thm:NV} applies to $u_S$. Let
$s_m:=\lfloor\eta^2m\rfloor$. For all sufficiently large $m$, one has
$s_m\ge C_1\log m$ and $\sqrt{s_m/m}\le\eta$. Consequently,
\[
\{S\notin\cG\}=\{\|u_S\|_\infty>\eta\}
\subseteq
\left\{\|u_S\|_\infty\ge\sqrt{\frac{s_m}{m}}\right\}.
\]
Applying Theorem~\ref{thm:NV} gives
\[
\Pp\{S\notin\cG\}\le Cm^2e^{-s_m/C}.
\]
For all sufficiently large $m$, $s_m\ge\eta^2m/2$, and hence
\[
\Pp\{S\notin\cG\}\le Cm^2e^{-cm}.
\]
Let
\[
N_{\mathrm{bad}}:=|\cS\setminus\cG|
=\sum_{S\in\cS}\1_{\{S\notin\cG\}}.
\]
By linearity of expectation,
\[
\Ee N_{\mathrm{bad}}
=\sum_{S\in\cS}\Pp\{S\notin\cG\}
\le LCm^2e^{-cm}.
\]
Since $|\cG|<L/2$ implies $N_{\mathrm{bad}}>L/2$, Markov's inequality yields
\[
\Pp\{|\cG|<L/2\}
\le\frac{2\Ee N_{\mathrm{bad}}}{L}
\le2Cm^2e^{-cm}.
\]
Absorbing the factor $2$ into $C$ proves the claim.
\end{proof}

We next control the geometry of pairs. The following entropy estimate shows that pairs of row sets
with very large overlap form an exponentially negligible fraction.

\begin{lemma}[Count high-overlap pairs]\label{lem:entropy}
Suppose $n_m/m\to\beta>1$ and let
\[
L_m:=\binom{n_m}{m-1}.
\]
Then
\begin{equation}\label{eq:L-entropy}
\frac1m\log L_m\longrightarrow h(\beta).
\end{equation}
Moreover, for every $\lambda<h(\beta)$ there are constants $\delta_0,\kappa>0$ such that, for all
sufficiently large $m$,
\begin{equation}\label{eq:overlap-entropy}
\frac1{L_m}
\sum_{1\le k\le\lfloor\delta_0m\rfloor}
\binom{m-1}{k}\binom{n_m-m+1}{k}
\le e^{-(\lambda+4\kappa)m}.
\end{equation}
\end{lemma}

\begin{proof}
Write $\beta_m:=n_m/m$, so $\beta_m\to\beta$. Stirling's formula gives
\[
\begin{aligned}
\frac1m\log L_m
&=\frac1m\bigl(\log(n_m!)-\log((m-1)!)-\log((n_m-m+1)!)\bigr)\\
&=\beta_m\log\beta_m-(\beta_m-1)\log(\beta_m-1)+o(1),
\end{aligned}
\]
which proves \eqref{eq:L-entropy}.

Fix $\lambda<h(\beta)$ and choose $\kappa>0$ so small that
$\lambda+7\kappa<h(\beta)$. By \eqref{eq:L-entropy}, for all sufficiently large $m$,
\[
L_m\ge e^{(h(\beta)-\kappa)m}.
\]
Since $n_m/m\to\beta$, there is $A>0$ such that $n_m-m+1\le Am$ for all sufficiently large $m$.
Using $\binom pk\le(ep/k)^k$, we obtain, for $1\le k\le\delta_0m$ and $x:=k/m$,
\[
\frac1m\log\left[\binom{m-1}{k}\binom{n_m-m+1}{k}\right]
\le x\log\left(\frac{Ae^2}{x^2}\right).
\]
The right-hand side tends to zero as $x\downarrow0$. Choose $\delta_0>0$ sufficiently small
that $\delta_0<\frac12\min\{1,\beta-1\}$ and
\[
x\log\left(\frac{Ae^2}{x^2}\right)\le\frac\kappa2
\qquad(0<x\le\delta_0).
\]
Then every summand in \eqref{eq:overlap-entropy} is at most $e^{\kappa m/2}$, and the sum contains
at most $m$ terms. Hence, for all sufficiently large $m$,
\[
\sum_{1\le k\le\lfloor\delta_0m\rfloor}
\binom{m-1}{k}\binom{n_m-m+1}{k}
\le me^{\kappa m/2}\le e^{\kappa m}.
\]
Dividing by the lower bound for $L_m$ gives
\[
\frac1{L_m}
\sum_{1\le k\le\lfloor\delta_0m\rfloor}
\binom{m-1}{k}\binom{n_m-m+1}{k}
\le e^{-(h(\beta)-2\kappa)m}.
\]
Since $h(\beta)-2\kappa>\lambda+4\kappa$, this proves \eqref{eq:overlap-entropy}.
\end{proof}
Lemma \ref{lem:entropy} shows that pairs $(S,T)$ with very large overlap, equivalently
with $|T\setminus S|$ small, form only a very small proportion of all
pairs. It remains to consider pairs for which $T\setminus S$ contains at
least $\delta_0m$ rows. Fix such a pair and condition on the rows indexed
by $S$. Then $u_S$ is fixed, while the rows indexed by $T\setminus S$
remain independent. If $u_T$ is close to $u_S$ up to a change of sign,
then, because the rows indexed by $T$ are orthogonal to $u_T$, the rows in
$T\setminus S$ must have small inner products with $u_S$, unless the matrix
formed by these rows has an unusually large operator norm. The following Lemma \ref{lem:operator-tail} rules
out the latter possibility with exponentially high probability.
Then Lemma \ref{lem:few-parallel} combines this operator-norm bound with the bounded-density
estimate for the projections of the rows indexed by $T\setminus S$ onto
$u_S$, showing that only an exponentially small proportion of the
remaining pairs can have close normal vectors.
\begin{lemma}[Operator-norm tail]\label{lem:operator-tail}
For every $H>0$ there is a constant $R=R(H,\xi)>0$ such that, for
every $1\le q\le m$ and every $q\times m$ matrix $B$ with i.i.d.
entries distributed as $\xi$,
\begin{equation}\label{eq:operator-tail}
 \Pp\{\|B\|_{\op}>R\sqrt m\}\le e^{-Hm}.
\end{equation}
\end{lemma}
\begin{proof}
The subgaussian operator-norm estimate
\cite[Theorem~4.4.5]{VershyninBook} gives a constant $C_\xi$ such that
\[
 \Pp\{\|B\|_{\op}>C_\xi(\sqrt q+\sqrt m+u)\}\le2e^{-u^2},
 \qquad u\ge0.
\]
Take $u=\sqrt{H+1}\sqrt m$. Since $q\le m$, one may choose
$R=C_\xi(2+\sqrt{H+1})$, and
$2e^{-(H+1)m}\le e^{-Hm}$ for $m\ge1$.
\end{proof}

\begin{lemma}[Few nearly parallel pairs]\label{lem:few-parallel}
Suppose $n_m/m\to\beta>1$. With the notation above, fix $\lambda<h(\beta)$. There are constants
$\delta_0,a,\kappa,c_*>0$ such that, with probability at least $1-e^{-c_*m}$ for all sufficiently
large $m$, the number of ordered distinct pairs $(S,T)\in\cS^2$ satisfying at least one of
\begin{equation}\label{eq:exceptional-pairs}
|T\setminus S|\le\delta_0m,
\qquad\text{or}\qquad
\dpr(u_S,u_T)<a,
\end{equation}
is at most
\[
2L_m^2e^{-(\lambda+4\kappa)m}.
\]
\end{lemma}

\begin{proof}
Fix $\lambda<h(\beta)$ and choose $\delta_0,\kappa>0$ from
Lemma~\ref{lem:entropy}. For each fixed $S\in\cS$, the number of
$T\in\cS$ with $|T\setminus S|=q$ is
$\binom{m-1}{q}\binom{n_m-m+1}{q}$. Hence the number of ordered
distinct pairs with $|T\setminus S|\le\delta_0m$ is at most
\[
 L_m^2e^{-(\lambda+4\kappa)m}.
\]
This estimate is deterministic.

Choose $H>\max\{0,\lambda+4\kappa\}+2$, and then choose
$R=R(H,\xi)$ from Lemma~\ref{lem:operator-tail}.
Fix an ordered pair $(S,T)$ with $q:=|T\setminus S|>\delta_0m$,
and condition on the rows indexed by $S$. Put $u=u_S$ and
$J=T\setminus S$, and let $B_J$ have the rows indexed by $J$.
These $q$ rows are independent of the conditioning and of one another.

If $\dpr(u_S,u_T)<a$, choose $v\in\{u_T,-u_T\}$ with
$\|u-v\|_2<a$. Since $B_Jv=0$, on
$\{\|B_J\|_{\op}\le R\sqrt m\}$ one has
\[
 \|B_Ju\|_2=\|B_J(u-v)\|_2<Ra\sqrt m.
\]
With $\tau_a:=Ra\sqrt{2/\delta_0}$, this implies
\[
 \#\{i\in J:|\langle R_i,u\rangle|\le\tau_a\}\ge q/2.
\]
Otherwise more than $q/2$ terms in $\|B_Ju\|_2^2$ would exceed
$\tau_a^2$, contradicting $q>\delta_0m$ and
$\|B_Ju\|_2^2<R^2a^2m$.

By Lemma~\ref{lem:marginal-density}, each of these small-ball
events has conditional probability at most
$p_a:=2\sqrt2K\tau_a$. Choose $a>0$ sufficiently small that
$p_a\le1$ and
$2\sqrt{p_a}\le\exp[-(H+1)/\delta_0]$.
Conditional independence and a union bound over subsets of rows give
\[
 \Pp\left\{
 \#\{i\in J:|\langle R_i,u\rangle|\le\tau_a\}\ge q/2
 \,\middle|\,(R_i)_{i\in S}\right\}
 \le2^q p_a^{q/2}
 \le e^{-(H+1)m}.
\]
Since $q\le m-1$, Lemma~\ref{lem:operator-tail} applies to $B_J$.
Thus, uniformly for these ordered pairs,
\[
 \Pp\{\dpr(u_S,u_T)<a\}\le2e^{-Hm}.
\]

Let $F_m$ count such low-overlap pairs with
$\dpr(u_S,u_T)<a$. Then
$\Ee F_m\le2L_m^2e^{-Hm}$, so Markov's inequality yields
\[
 \Pp\{F_m>L_m^2e^{-(\lambda+4\kappa)m}\}
 \le2e^{-(H-\lambda-4\kappa)m}
 \le e^{-c_*m}
\]
for all sufficiently large $m$, with
$c_*=(H-\lambda-4\kappa)/2>0$. Combining this with the
deterministic high-overlap count proves the lemma.
\end{proof}

The upper estimate also relies on estimation of some small ball probabilities. The following two corollaries give the results in the subgaussian setting.
\begin{corollary}[One-slab estimate for marked rows]\label{cor:one-slab}
Let $X_1,\ldots,X_k\in\R^m$ be independent rows with i.i.d.
coordinates distributed as $\xi$, and let the unit vectors below
be deterministic. Define the following small ball probability
$$p_{u,k}(t):=\Pp\left\{\sum_{j=1}^k|\langle X_j,u\rangle|^2\le t^2\right\}$$
There are constants $c_1,C_1>0$ such that for every $k\ge1$, every $0<t\le1$, and every unit vector $u$ with $\|u\|_\infty\le\eta$,
\begin{equation}\label{eq:one-normal-marked}
 c_1^k v_kt^k\le p_{u,k}(t)\le C_1^k v_kt^k,
\end{equation}
where $v_q$ is the volume of the Euclidean unit ball in $\mathbb R^q$.
\end{corollary}
\begin{proof}
    By Lemma \ref{lem:local-lower} and independence, the joint density of $(\langle X_1,u\rangle,\ldots,\langle X_k,u\rangle)$ is at least $c_0^k$ on the unit ball in $\mathbb R^k$, integration over the ball of radius $t$ gives the lower bound. The upper bound follows similarly from \eqref{eq:one-dimensional-density}.
\end{proof}
\begin{corollary}[Two-slab estimate for marked rows]\label{cor:two-slab}
Let $X_1,\ldots,X_k\in\R^m$ be independent rows with i.i.d.
coordinates distributed as $\xi$, and let the unit vectors below
be deterministic. Define
$$q_{u,v,k}(t):=\Pp\left\{\sum_{j=1}^k|\langle X_j,u\rangle|^2\le t^2,\sum_{j=1}^k|\langle X_j,v\rangle|^2\le t^2\right\}$$
For arbitrary unit vectors $u,v$ with $\|u\|_{\infty},\|v\|_{\infty}\le\eta$, and any $0<t\leq1$,
\begin{equation}\label{eq:two-normal-crude}
 q_{u,v,k}(t)\le C_1^k v_kt^k.
\end{equation}
Furthermore, assume the above $u,v$ also satisfies $\dpr(u,v)\ge a>0$, then
\begin{equation}\label{eq:two-normal-marked}
 q_{u,v,k}(t)\le C_a^k v_k^2t^{2k}.
\end{equation}
where $C_a$ is a constant which only depends on $a,C_{\mathrm{RV}},K$.
\end{corollary}
\begin{proof}
The crude bound $\eqref{eq:two-normal-crude}$ follows by discarding one event and the right-hand side of $\eqref{eq:one-normal-marked}$. For the sharper bound, note that if $\rho=|\langle u,v\rangle|$, then
$$d_{\mathrm{pr}}^2(u,v)=2(1-\rho)$$
Hence if $\dpr(u,v)\ge a$, then
$$1-\rho^2=(1-\rho)(1+\rho)\geq 1-\rho\geq a^2/2$$
Thus, by \eqref{eq:two-dimensional-density}, we have for every $j$,
$$\|\dens(\langle X_j,u\rangle,\langle X_j,v\rangle)\|_\infty
\le \frac{C_{\mathrm{RV}}^2K^2}{\sqrt{1-\rho^2}}\leq \frac{\sqrt{2}C^2_{\mathrm{RV}}K^2}{a}:=C_a
$$
Thus similarly a integration over ball of radius $t$ gives the desired result.
\end{proof}
\section{Sequential marked rows and the upper estimate}
\label{sec:sequential}

We now prove the upper estimate. As in the usual marked-row construction, an $(m-1)$-tuple of base rows determines a candidate hyperplane, and one seeks marked rows that are almost orthogonal to its normal. The difference is that the $k$ marked constraints are imposed one at a time. A one-row second-moment estimate has a success probability bounded away from zero, and the candidate family is then restricted successively using $k$ independent groups of reserved rows.

\subsection{Proof map.}
Reserve $kb\asymp m$ rows and split them into $k$ groups of
size $b\asymp m/k$. The remaining rows generate exponentially many
candidate hyperplanes. After selecting $j$ marked rows, retain those
candidates whose normal has scalar product at most $s=t/\sqrt{k}$
with each selected row. Lemma~\ref{lem:one-row-stage} shows that a
fresh row has a fixed positive probability of retaining a proportion
at least $c_1s$ of the current candidates while leaving at most
$A_0s$ times the current number of exceptional pairs.
Lemma~\ref{lem:iteration-stability} verifies that the required bounds
persist through all $k$ stages. Independence within each reserved
group then gives a suitable row with probability at least $1-e^{-cb}$.

\subsection{Sequential filtering of candidate hyperplanes}

Let $b_1,\ldots,b_n\in\R^m$ be the base rows, and let
$x_1,\ldots,x_k\in\R^m$ be distinct marked rows chosen outside this
base family. In the probabilistic construction below, these marked
rows will be selected from $k$ disjoint reserved groups; the other
reserved rows are not part of the base family. Set
\[
 \cS:=\{S\subset[n]:|S|=m-1\},
\]
and let $u_S$ be a unit normal to the base rows indexed by $S$.
For $\mathcal I\subset\cS$, $x\in\R^m$, and $s>0$, define
\begin{equation}\label{eq:one-row-restriction}
 \mathcal I(x;s):=\{S\in\mathcal I:|\langle x,u_S\rangle|\le s\}.
\end{equation}
Starting from $\mathcal I_0:=\mathcal I$, define
\begin{equation}\label{eq:candidate-chain}
 \mathcal I_j:=\mathcal I_{j-1}(x_j;s),\qquad 1\le j\le k.
\end{equation}

Hence at the last iteration, if $\mathcal I_k$ is nonempty, there is a set $S^*$, such that for every $j=1,2,\ldots,k$, $|\langle x_j,u_{S^*}\rangle|\leq s$.
\begin{lemma}[Sequential marked rows imply a bad submatrix]\label{lem:sequential-deterministic}

Let $x_1,\ldots,x_k$ be the $k$ marked rows, let $t>0$, and set $s=t/\sqrt{k}$. If the family $\mathcal I_k$ in \eqref{eq:candidate-chain} is nonempty, then the union of the $m-1$ base rows indexed by $S^*$ and the $k$ marked rows consists of an $(m+k-1)\times m$ row submatrix $B$ satisfying
$\sigma_m(B)\leq t$
and thus
\[
 M_{m,k-1}(A_m)\le t.
\]
\end{lemma}

\begin{proof}
Since the base rows are orthogonal to $u_{S^*}$ and
\[
 |\langle x_j,u_{S^*}\rangle|\le s,
 \qquad 1\le j\le k.
\]
We have
\[
 \|Bu_{S^*}\|_2^2
 =\sum_{i\in S^*}|\langle b_i,u_{S^*}\rangle|^2+\sum_{j=1}^k|\langle x_j,u_{S^*}\rangle|^2
 \le ks^2=t^2.
\]
The variational characterization of the least singular value yields $\sigma_m(B)\le t$ and thus $M_{m,k-1}(A_m)\le t$.
\end{proof}
\begin{lemma}[Paley--Zygmund inequality]\label{lem:PZ}
For every
nonnegative random variable $Z$ with $0<\Ee Z^2<\infty$ and every
$\rho\in(0,1)$,
\[
\mathbb P\{Z\ge\rho\,\mathbb EZ\}
\ge
(1-\rho)^2
\frac{(\mathbb EZ)^2}{\mathbb EZ^2}.
\]
\end{lemma}
\begin{lemma}[One marked-row stage]\label{lem:one-row-stage}
Let $\mathcal F\subset\cS$ be finite, with
$\|u_S\|_\infty\le\eta$ for every $S\in\mathcal F$. Let
\[
 \mathcal E\subset\{(S,T)\in\mathcal F^2:S\ne T\}
\]
contain every distinct pair for which $\dpr(u_S,u_T)<a$.
Write $M=|\mathcal F|$, and assume
\begin{equation}\label{eq:thinning}
 sM\ge1,\qquad |\mathcal E|\le sM^2,\qquad 0<s\le1.
\end{equation}
Let $X$ have i.i.d. coordinates distributed as $\xi$, independently
of the family and its normals. Define
\[
 Z(X):=\sum_{S\in\mathcal F}
       \1_{\{|\langle X,u_S\rangle|\le s\}},
\qquad
 W(X):=\sum_{(S,T)\in\mathcal E}
       \1_{\{|\langle X,u_S\rangle|\le s,\,
              |\langle X,u_T\rangle|\le s\}}.
\]
There are constants $c_1\in(0,1)$, $p_0>0$, and $A_0>1$, depending
only on $a$ and the law of $\xi$, such that, conditionally on the
family and its normals,
\[
 \Pp\{Z(X)\ge c_1sM,\ W(X)\le A_0s|\mathcal E|\}\ge p_0.
\]
\end{lemma}
\begin{proof}
Condition on the family and its normals. The one-row cases of
Corollaries~\ref{cor:one-slab} and \ref{cor:two-slab} give
\[
 \Ee Z\ge csM,\qquad
 \Ee Z^2\le CsM+C_as^2M^2+Cs|\mathcal E|
 \le C's^2M^2.
\]
The three terms in the second-moment estimate correspond to diagonal
pairs, distinct pairs outside $\mathcal E$, and pairs in $\mathcal E$.
The final inequality uses \eqref{eq:thinning}.
Paley--Zygmund therefore gives
\[
 \Pp\{Z\ge\tfrac12\Ee Z\}\ge\frac{c^2}{4C'}=:p_1>0.
\]
Choose $c_1=\min\{c/2,1/2\}$. Then
$\Pp\{Z\ge c_1sM\}\ge p_1$.

The crude one-vector bound gives $\Ee W\le Cs|\mathcal E|$.
If $\mathcal E=\varnothing$, then $W=0$ identically.
Otherwise Markov's inequality, with $A_0>1$ chosen so that
$C/A_0\le p_1/2$, gives
\[
 \Pp\{W>A_0s|\mathcal E|\}\le p_1/2.
\]
In either case the simultaneous event has probability at least
$p_0:=p_1/2$.
\end{proof}

\subsection{Iteration estimates}

Fix $\eps\in(0,h(\gamma))$ and, throughout the remainder of the section, write
\begin{equation}\label{eq:upper-parameters}
 k:=k_m,
 \qquad
 \lambda:=h(\gamma)-\eps,
 \qquad
 t:=\exp\left(-\lambda\frac{m}{k}\right),
 \qquad
 s:=\frac{t}{\sqrt{k}}.
\end{equation}
Since $h(\cdot)$ is increasing, we can choose a small enough $\delta\in(0,\gamma-1)$ such that, with $\beta:=\gamma-\delta$,
\begin{equation}\label{eq:entropy-gap}
 \Delta:=h(\beta)-\lambda=h(\gamma-\delta)-h(\gamma)+\varepsilon>0.
\end{equation}
Now set
$
b:=\left\lfloor\frac{\delta m}{k}\right\rfloor.
$
and choose $kb$ rows and split them into $k$ disjoint groups $\{\mathcal R_1,\ldots,\mathcal{R}_k\}$, where
\[
 \mathcal R_j=\{Y_{j,1},\ldots,Y_{j,b}\},
 \quad 1\le j\le k.
\]
The remaining
$
 n:=N_m-kb
$
rows are the base rows. Since $kb=\delta m+O(k)$ and $k=o(m)$, one has $n/m\to\beta>1$.

Now let
\[
 \cS=\{S\subset[n]:|S|=m-1\},
 \qquad
 L:=|\cS|=\binom{n}{m-1}.
\]
Since $\lambda<h(\beta)$, Lemma~\ref{lem:many-delocalized} and Lemma~\ref{lem:few-parallel} implies there are constants $a,\kappa,c>0$ and a event $\Omega_{\mathrm{base}}$ such that
\begin{equation}\label{eq:base-event-probability}
 \Pp(\Omega_{\mathrm{base}}^c)\le Ce^{-cm},
\end{equation}
and, on $\Omega_{\mathrm{base}}$, there is a subset $\cG\subset\cS$ and a set of ordered distinct pairs
$$\mathcal{E}=\{(S,T)\in\mathcal{G}^2:S\neq T,|T\setminus S|\leq\delta_0m\text{ or } d_{\mathrm{pr}}(u_S,u_T)<a\}$$
satisfying
\begin{equation}\label{eq:base-event-bounds}
 |\cG|\ge\frac L2,
 \qquad
 |\cE|\le2L^2e^{-(\lambda+4\kappa)m},
\end{equation}
while every distinct pair outside $\cE$ has projective distance $d_{\mathrm{pr}}(u_S,u_T)$ at least $a$.

Starting from $\mathcal I_0:=\cG$, suppose marked row $x_j$ have been chosen successively and define
\[
 \mathcal I_j:=\mathcal I_{j-1}(x_j;s),
 \qquad
 M_j:=|\mathcal I_j|,
 \qquad
 Q_j:=|\cE\cap(\mathcal I_j\times\mathcal I_j)|.
\]
Here $M_j$ measures the number of candidate hyperplanes that survive after the first $j$-th round. And $Q_j$ measures the number of surviving pairs whose normal vectors may be too strongly correlated. Controlling $Q_j$ ensures that such pairs remain sufficiently sparse for the second-moment argument at the next stage.
\begin{lemma}[Stability of the iteration]\label{lem:iteration-stability}
Assume $k\log m=o(m)$ and that $\Omega_{\mathrm{base}}$ holds.
Use the constants $c_1,A_0$ from Lemma~\ref{lem:one-row-stage}.
Suppose that the first $j$ selected marked rows satisfy
\begin{equation}\label{eq:successful-thinning}
 M_i\ge c_1sM_{i-1},\qquad Q_i\le A_0sQ_{i-1},
 \qquad 1\le i\le j.
\end{equation}
Then
\begin{equation}\label{eq:iteration-bounds}
 M_j\ge\frac L2(c_1s)^j,\qquad
 Q_j\le2L^2e^{-(\lambda+4\kappa)m}(A_0s)^j.
\end{equation}
Uniformly for $0\le j<k$, these bounds imply
$sM_j\ge1$ and $Q_j\le sM_j^2$ for all sufficiently large $m$.
If all $k$ stages succeed, then $M_k\ge1$.
\end{lemma}
\begin{proof}
Iteration of \eqref{eq:successful-thinning}, starting from
\eqref{eq:base-event-bounds}, gives \eqref{eq:iteration-bounds}.
Also,
\begin{equation}\label{eq:base-family-entropy}
 \log L=h(\beta)m+o(m).
\end{equation}
Recall that $\Delta=h(\beta)-\lambda>0$,
$\log s=-\lambda m/k-\tfrac12\log k$, and $0<c_1<1$.
For $0\le j<k$,
\begin{align}
 \log(sM_j)
 &\ge\log(L/2)+(j+1)\log s+j\log c_1 \notag\\
 &\ge\Delta m-\frac k2\log k-O(k)+o(m)
 \longrightarrow+\infty .
 \label{eq:candidate-budget}
\end{align}
The lower bound is uniform in $j$. Similarly, since $A_0>1$,
\begin{align}
 \frac{Q_j}{sM_j^2}
 &\le8e^{-(\lambda+4\kappa)m}
       \left(\frac{A_0}{c_1^2}\right)^j s^{-j-1} \notag\\
 &\le\exp\left\{-4\kappa m+\frac k2\log k+O(k)\right\}
 \longrightarrow0,
 \label{eq:pair-budget}
\end{align}
uniformly for $j<k$, using
$s^{-k}=e^{\lambda m}k^{k/2}$.
These limits follow from $k\log m=o(m)$ and imply the stage
hypotheses. Finally,
\[
 \log M_k
 \ge\log(L/2)+k\log s+k\log c_1
 \ge\Delta m-\frac k2\log k-O(k)+o(m)
 \longrightarrow+\infty.
\]
Thus $M_k\ge1$ when all stages succeed.
\end{proof}

\subsection{Completion of the upper estimate}

\begin{proposition}[Upper estimate]\label{prop:sequential-upper}
Under the assumptions of Theorem \ref{thm:moderate}, for every fixed $\eps>0$, there are constants $c,C>0$ such that, with
\[
 t_m:=\exp\left[-(h(\gamma)-\eps)\frac{m}{k_m}\right],
\]
one has, for all sufficiently large $m$,
\begin{equation}\label{eq:sequential-upper}
 \Pp\{M_{m,r_m}(A_m)>t_m\}
 \le Ce^{-cm}+k_m\exp\left(-c\frac{m}{k_m}\right).
\end{equation}
In particular, the probability on the left tends to zero.
\end{proposition}

\begin{proof}
Without loss of generality, we assume $\eps\in(0,h(\gamma))$. Use the notation and row decomposition introduced above, and condition first on the base rows. On $\Omega_{\mathrm{base}}$, set $\mathcal I_0=\cG$. We choose one row from each group inductively.

Suppose that rows $x_1,\ldots,x_{j-1}$ have already been chosen and that the first $j-1$ thinnings satisfy \eqref{eq:successful-thinning}. The family $\mathcal I_{j-1}$ is determined by the base rows and the first $j-1$ groups, whereas the rows in the $j$th group $\mathcal R_j$ are independent of this information. By Lemma \ref{lem:iteration-stability}, $\mathcal I_{j-1}$ satisfies \eqref{eq:thinning}. Lemma \ref{lem:one-row-stage} therefore shows that each row in $\mathcal R_j$, conditionally on $\mathcal{F}_{j-1}=\sigma(\text{base rows},\mathcal{R}_1,\ldots,\mathcal{R}_{j-1})$, satisfies \eqref{eq:successful-thinning} with probability at least $p_0$. Since the $b$ rows in this group are conditionally independent,
\begin{equation}\label{eq:one-group-failure}
 \Pp\{\mathcal R_j\text{ contains no suitable row}\mid\mathcal{F}_{j-1}\}
 \le(1-p_0)^b\le e^{-p_0b}.
\end{equation}
Whenever a suitable row exists, choose the first one in the fixed ordering of $\mathcal R_j$ and continue the construction.

A union bound over the $k$ groups gives
\begin{equation}\label{eq:all-groups-failure}
 \Pp\{\text{the construction stops before round }k\mid\Omega_{\mathrm{base}}\}
 \le ke^{-p_0b}.
\end{equation}
If the construction reaches round $k$, Lemma~\ref{lem:iteration-stability} gives $\mathcal I_k\ne\varnothing$, and Lemma~\ref{lem:sequential-deterministic} then yields $M_{m,r_m}(A_m)\le t_m$.

Because $b=\lfloor\delta m/k\rfloor$ and $m/k\to\infty$, one has $b\ge\delta m/(2k)$ for all sufficiently large $m$. Combining \eqref{eq:base-event-probability} and \eqref{eq:all-groups-failure} therefore gives
\[
 \Pp\{M_{m,r_m}(A_m)>t_m\}
 \le Ce^{-cm}+k\exp\left(-c\frac{m}{k}\right).
\]
Finally, $k\log m=o(m)$ implies $\log k=o(m/k)$, so the second term tends to zero.
\end{proof}

\begin{proof}[Proof of Theorem~\ref{thm:moderate}]
Proposition~\ref{prop:lower-tail} gives the lower deviation
\eqref{eq:lower-tail-event}, and
Proposition~\ref{prop:sequential-upper} gives the upper deviation
\eqref{eq:sequential-upper}. For every fixed
$\varepsilon\in(0,h(\gamma))$, their combination yields
\[
 \Pp\left\{
 \left|\frac{k_m}{m}\log M_{m,r_m}(A_m)+h(\gamma)\right|
 >\varepsilon\right\}\longrightarrow0.
\]
This proves \eqref{eq:th1.1}.
\end{proof}

\begin{proof}[Proof of Theorem~\ref{thm:main}]
Take $k_m=r+1$ fixed in Theorem~\ref{thm:moderate} and divide the
limit by $r+1$. This gives \eqref{eq:fixed}. For each fixed
deviation tolerance, \eqref{eq:lower-tail-event} and
\eqref{eq:sequential-upper} also give exponentially small failure
probabilities.
\end{proof}

\section{Phase retrieval}

We now apply the row-submatrix estimates to the stability of real phase retrieval. As explained in the Introduction, the loss of sign information makes this stability depend on the least singular values of complementary row blocks. For a full-spark matrix with $N=2m-1+r$ rows, the obstruction $\omega(A)$ introduced below is exactly the row minimum $M_{m,r}(A)$. At the critical count $r=0$, this also equals the optimal lower Lipschitz constant $L_A$. When $r>0$, partitions with two full-rank blocks require an additional estimate. We combine the row-minimum theorem with small-ball estimates for these complementary blocks to determine the logarithmic decay rate of $L_A$ and the corresponding growth rate of the phase-retrieval condition number.\\

\smallskip
\noindent
Recall that, for $A\in\R^{N\times m}$, the real amplitude map is
\[
 \Phi_A:\R^m/\{\pm1\}\longrightarrow\R^N,\qquad
 \Phi_A(x)=|Ax|,
\]
where the absolute value is taken coordinatewise. Equip the domain
with the quotient distance
\[
 d_\pm(x,y):=\min\{\|x-y\|_2,\|x+y\|_2\}.
\]
The optimal lower and upper Lipschitz constants are
\begin{equation}\label{eq:def-LU}
 L_A:=\inf_{x\ne\pm y}
       \frac{\|\Phi_A(x)-\Phi_A(y)\|_2}{d_\pm(x,y)},
 \qquad
 U_A:=\sup_{x\ne\pm y}
       \frac{\|\Phi_A(x)-\Phi_A(y)\|_2}{d_\pm(x,y)}.
\end{equation}
Define also the rank-deficient-side obstruction
\begin{equation}\label{eq:def-omega}
 \omega(A):=
 \min_{\substack{J\subset[N]\\ \rank(A_J)<m}}
       \sigma_m(A_{J^c}).
\end{equation}
Here $J^c=[N]\setminus J$, and the variational convention for
$\sigma_m$ in Section~1 includes empty and rank-deficient blocks.
When $L_A>0$, the amplitude-map condition number is
\begin{equation}\label{eq:def-phase-condition}
 \kappa_{\mathrm{ph}}(A):=\frac{U_A}{L_A}
 =\frac{\|A\|_{\op}}{L_A}.
\end{equation}

\begin{proposition}[Partition formula]\label{prop:phase-partition}
For every real $N\times m$ matrix $A$,
\begin{equation}\label{eq:phase-partition}
 L_A=\min_{J\subset[N]}
       \left(\sigma_m(A_J)^2+\sigma_m(A_{J^c})^2\right)^{1/2},
 \qquad U_A=\|A\|_{\op}.
\end{equation}
Consequently $L_A\le\omega(A)$. If $N=2m-1$, then
$L_A=\omega(A)$.
\end{proposition}
\begin{proof}
The first identity is the real amplitude-map formula of
\cite[equations~(4.3)--(4.5) and Theorem~4.3]{BalanWang};
we give a direct verification in row-matrix notation.
For real $a,b$,
\[
 \bigl||a|-|b|\bigr|^2=\min\{|a-b|^2,|a+b|^2\}.
\]
Writing $u=x-y$ and $v=x+y$ therefore gives
\[
 \|\Phi_A(x)-\Phi_A(y)\|_2^2
 =\sum_{i=1}^N
   \min\{|\langle a_i,u\rangle|^2,|\langle a_i,v\rangle|^2\}
 =\min_{J\subset[N]}
   \bigl(\|A_Ju\|_2^2+\|A_{J^c}v\|_2^2\bigr),
\]
where $a_i$ is the $i$th row of $A$.
Let $\Lambda$ denote the minimum of
$\sigma_m(A_J)^2+\sigma_m(A_{J^c})^2$ over $J$.
Every term in the last minimum is at least
$\Lambda\min\{\|u\|_2^2,\|v\|_2^2\}$, so $L_A^2\ge\Lambda$.

Conversely, fix $J$ and choose unit vectors $u,v$ attaining
$\sigma_m(A_J)$ and $\sigma_m(A_{J^c})$, respectively.
Set $x=(u+v)/2$ and $y=(v-u)/2$.
Then $d_\pm(x,y)=1$ and $x\ne\pm y$, so the same identity gives
$L_A^2\le\sigma_m(A_J)^2+\sigma_m(A_{J^c})^2$.
Minimizing over $J$ proves the first formula.
The bound $U_A\le\|A\|_{\op}$ follows from the reverse triangle
inequality, applied with both $y$ and $-y$; equality follows by
taking $y=0$ and a unit vector $x$ attaining $\|A\|_{\op}$.

Restricting the partition minimum to $\rank(A_J)<m$ gives
$L_A\le\omega(A)$. When $N=2m-1$, every partition has a side
with at most $m-1$ rows. That side has rank less than $m$, so,
after exchanging the sides if necessary, every value in the
partition minimum occurs in \eqref{eq:def-omega}. Hence equality
holds.
\end{proof}

For $N\ge2m$, a partition can have two full-rank sides, so equality
with $\omega(A)$ no longer follows. The next theorem identifies the
rate of $\omega(A)$; Corollary~\ref{cor:phase-lower-Lipschitz} then
establishes the matching rate for $L_A$ in the moderate-surplus regime.

\begin{definition}[Full spark]
A real $N\times m$ matrix is full spark if every $m$ rows are linearly independent.
\end{definition}

\begin{lemma}[Reduction of $\omega$ to the redundant row minimum]
\label{lem:omega-reduction}
Let $r\ge0$ be an integer, and let
$A\in\R^{(2m-1+r)\times m}$ be full spark. Then
\begin{equation}\label{eq:omega-reduction}
 \omega(A)=M_{m,r}(A).
\end{equation}
Moreover, $L_A>0$.
\end{lemma}
\begin{proof}
Full spark implies $\rank(A_J)<m$ exactly when $|J|\le m-1$.
If $J_1\subset J_2$, deleting rows gives
$\sigma_m(A_{J_2^c})\le\sigma_m(A_{J_1^c})$.
Thus the minimum in \eqref{eq:def-omega} is attained among the
sets with $|J|=m-1$, whose complements have $m+r$ rows.
This proves \eqref{eq:omega-reduction}.
Every partition of $2m-1+r$ rows has at least one side with
at least $m$ rows; full spark makes its least singular value
positive. The finite minimum in \eqref{eq:phase-partition}
is therefore positive.
\end{proof}

\begin{theorem}[Phase-retrieval obstruction rate]\label{thm:phase}
Let $r_m\ge0$ be integers, put $k_m=r_m+1$ and
$N_m=2m-1+r_m$, and assume $k_m\log m=o(m)$.
Let $A_m\in\R^{N_m\times m}$ have i.i.d. entries distributed as
$\xi$ under the standing assumptions. Then $A_m$ is full spark
almost surely, and
\begin{equation}\label{eq:phase-omega-rate}
 \frac{k_m}{m}\log\omega(A_m)
 \xrightarrow{\Pp}-\log4.
\end{equation}
\end{theorem}
\begin{proof}
Absolute continuity gives full spark almost surely.
Lemma~\ref{lem:omega-reduction} yields
$\omega(A_m)=M_{m,r_m}(A_m)$.
Since $N_m=2m+k_m-2$ and $k_m=o(m)$, one has $N_m/m\to2$.
Theorem~\ref{thm:moderate}, together with $h(2)=\log4$, proves
\eqref{eq:phase-omega-rate}.
\end{proof}

\begin{corollary}[Optimal lower Lipschitz rate with redundant rows]
\label{cor:phase-lower-Lipschitz}
Under the assumptions of Theorem~\ref{thm:phase},
\begin{equation}\label{eq:phase-L-rate}
 \frac{k_m}{m}\log L_{A_m}\xrightarrow{\Pp}-\log4.
\end{equation}
More precisely, with $D$ from
Lemma~\ref{lem:rectangular-small-ball}, for every $t>0$,
\begin{equation}\label{eq:phase-L-small-ball}
 \Pp\{L_{A_m}\le t\}
 \le2^{N_m}m^2(D\sqrt m\,t)^{k_m}.
\end{equation}
In particular, for every fixed $\varepsilon>0$,
\begin{equation}\label{eq:phase-L-lower-tail}
 \Pp\left\{L_{A_m}\le
   \exp\left[-(\log4+\varepsilon)\frac{m}{k_m}\right]\right\}
 \le\exp\{-\varepsilon m+O(k_m\log m)\}\longrightarrow0.
\end{equation}
For fixed $r_m=r$, this gives
$m^{-1}\log L_{A_m}\xrightarrow{\Pp}-\log4/(r+1)$.
\end{corollary}
\begin{proof}
Fix $m$ and abbreviate $N=N_m$, $k=k_m$, and $A=A_m$.
For a deterministic $J\subset[N]$, let
\[
 F_J(t):=\{\sigma_m(A_J)^2+\sigma_m(A_{J^c})^2\le t^2\}.
\]
By \eqref{eq:phase-partition},
$\{L_A\le t\}=\bigcup_{J\subset[N]}F_J(t)$.
We bound the probability for each fixed partition.

Suppose first that $|J|\ge m$ and $|J^c|\ge m$.
Set
\[
 d_1:=|J|-m+1,\qquad d_2:=|J^c|-m+1.
\]
Then $d_1,d_2\ge1$ and, crucially,
\[
 d_1+d_2=N-2m+2=k.
\]
The event $F_J(t)$ requires both $\sigma_m(A_J)\le t$ and
$\sigma_m(A_{J^c})\le t$. The two row blocks are independent
because $J$ is deterministic and the entries of $A$ are independent.
Applying Lemma~\ref{lem:rectangular-small-ball} separately gives
\begin{align*}
 \Pp(F_J(t))
 &\le\Pp\{\sigma_m(A_J)\le t\}
       \Pp\{\sigma_m(A_{J^c})\le t\}\\
 &\le m^2(D\sqrt m\,t)^{d_1+d_2}
 =m^2(D\sqrt m\,t)^k.
\end{align*}

Suppose instead that $|J|<m$. Then $\sigma_m(A_J)=0$ and
$|J^c|\ge m+k-1$. Choose a deterministic subset
$T_J\subset J^c$ of size $m+k-1$.
Adding rows cannot decrease the least singular value, so
$\sigma_m(A_{J^c})\ge\sigma_m(A_{T_J})$. Hence
\[
 \Pp(F_J(t))
 \le\Pp\{\sigma_m(A_{T_J})\le t\}
 \le m(D\sqrt m\,t)^k
 \le m^2(D\sqrt m\,t)^k.
\]
The case $|J^c|<m$ is identical after interchanging the sides.
These cases exhaust all partitions since $N\ge2m-1$.
In particular, when $k=1$ there are no partitions with both
sides of size at least $m$.

A union bound over the $2^N$ subsets $J$ proves
\eqref{eq:phase-L-small-ball}; no independence between different
partitions is needed.
For $t=\exp[-(\log4+\varepsilon)m/k]$, the logarithm of its
right-hand side is
\begin{align*}
 N\log2+2\log m+k\log D+\frac k2\log m
 -(\log4+\varepsilon)m
 &=-\varepsilon m+(k-2)\log2+2\log m\\
 &\quad+k\log D+\frac k2\log m\\
 &=-\varepsilon m+O(k\log m).
\end{align*}
The assumption $k\log m=o(m)$ proves
\eqref{eq:phase-L-lower-tail}.

For the opposite deviation, use $L_A\le\omega(A)$ and
Theorem~\ref{thm:phase}. More quantitatively, full spark and
Proposition~\ref{prop:sequential-upper}, with $\gamma=2$, give,
for each fixed $\varepsilon\in(0,\log4)$,
\begin{equation}\label{eq:phase-L-upper-tail}
 \Pp\left\{L_{A_m}>
 \exp\left[-(\log4-\varepsilon)\frac{m}{k_m}\right]\right\}
 \le Ce^{-cm}+k_m\exp\left(-c\frac{m}{k_m}\right)
 \longrightarrow0.
\end{equation}
Combining the two deviations proves \eqref{eq:phase-L-rate}.
\end{proof}

\begin{remark}[Scope of the rate conclusion]
Theorem~\ref{thm:phase} and
Corollary~\ref{cor:phase-lower-Lipschitz} identify the same
logarithmic rate for $\omega(A_m)$ and $L_{A_m}$.
They do not assert equality of these quantities for redundant
matrices, convergence of their ratio to one, or a sharp
multiplicative prefactor. The conclusion is restricted to
$k_m\log m=o(m)$.
\end{remark}

\begin{corollary}[Almost-sure logarithmic laws]\label{cor:almost-sure}
Assume the matrices $A_m$ in Theorem~\ref{thm:moderate} are defined
on a common probability space. No independence across dimensions
is required. Then
\[
 \frac{k_m}{m}\log M_{m,r_m}(A_m)\longrightarrow-h(\gamma)
 \qquad\text{almost surely}.
\]
In particular, for fixed $r_m=r$, the limit in
Theorem~\ref{thm:main} holds almost surely.
Under the phase-retrieval specialization $N_m=2m-1+r_m$,
\[
 \frac{k_m}{m}\log\omega(A_m)\longrightarrow-\log4,
 \qquad
 \frac{k_m}{m}\log L_{A_m}\longrightarrow-\log4
 \qquad\text{almost surely}.
\]
These statements apply to each prescribed sequence satisfying the
stated assumptions.
\end{corollary}
\begin{proof}
Fix $\varepsilon\in(0,h(\gamma))$. By
\eqref{eq:lower-tail-event}, the lower-deviation probability is at
most $\exp\{-\varepsilon m+o(m)\}$ and is therefore summable.
By \eqref{eq:sequential-upper}, the upper-deviation probability is
at most
\[
 Ce^{-cm}+k_m\exp(-cm/k_m),
\]
where $c,C>0$ may depend on the fixed $\varepsilon$ and $\gamma$,
but not on $m$. Since $k_m\log m=o(m)$, eventually $k_m\le m$
and $cm/k_m\ge4\log m$. Hence
\[
 k_m\exp(-cm/k_m)\le m^{-3}
\]
for all sufficiently large $m$, so the upper-deviation probabilities
are summable as well. The first Borel--Cantelli lemma, applied for
a countable sequence of positive $\varepsilon$ tending to zero,
proves the row-minimum limit. Taking $k_m=r+1$ gives the fixed-surplus
assertion. Absolute continuity gives the required full-rank events
simultaneously for all $m$ by a countable intersection.

For $\omega(A_m)$, use \eqref{eq:omega-reduction} and $h(2)=\log4$.
For $L_{A_m}$, the lower-deviation bound
\eqref{eq:phase-L-lower-tail} is summable because $k_m\log m=o(m)$,
and the upper-deviation bound \eqref{eq:phase-L-upper-tail} is
summable by the same argument. Applying Borel--Cantelli again
proves the remaining limit. None of these applications requires
independence across dimensions.
\end{proof}

\begin{corollary}[Phase-retrieval condition number]
\label{cor:phase-condition-rate}
Under the assumptions of Theorem~\ref{thm:phase},
\begin{equation}\label{eq:phase-condition-rate}
 \frac{k_m}{m}\log\kappa_{\mathrm{ph}}(A_m)
 \xrightarrow{\Pp}\log4.
\end{equation}
If the matrices are defined on a common probability space, this
convergence holds almost surely, without requiring independence
across dimensions. Moreover,
$\kappa_{\mathrm{ph}}(A_m/\sqrt m)=\kappa_{\mathrm{ph}}(A_m)$.
\end{corollary}
\begin{proof}
Since $N_m/m\to2$, the subgaussian operator-norm estimate
\cite[Theorem~4.4.5]{VershyninBook}, with deviation parameter
$u=\sqrt m$, gives constants $c,C>0$ such that
\[
 \Pp\{\|A_m\|_{\op}>C\sqrt m\}\le2e^{-cm}
\]
for all sufficiently large $m$.
For a lower bound, $\|A_m\|_{\op}\ge\|A_me_1\|_2$.
Because $\Ee\xi^2=1$, there is a fixed $\theta>0$ such that
\[
 \Ee e^{-\theta\xi^2}\le e^{-3\theta/4}.
\]
Indeed, both sides equal one at $\theta=0$, and their right
derivatives there are $-1$ and $-3/4$, respectively.
Independence of the entries in the first column and Markov's
inequality give
\[
 \Pp\{\|A_me_1\|_2^2\le N_m/2\}
 \le e^{\theta N_m/2}\bigl(\Ee e^{-\theta\xi^2}\bigr)^{N_m}
 \le e^{-\theta N_m/4}.
\]
Thus $\log\|A_m\|_{\op}=\tfrac12\log m+O_{\Pp}(1)$, and
\[
 \frac{k_m}{m}\log\|A_m\|_{\op}\xrightarrow{\Pp}0
\]
by $k_m\log m=o(m)$. Combining this with
\eqref{eq:def-phase-condition} and
Corollary~\ref{cor:phase-lower-Lipschitz} proves
\eqref{eq:phase-condition-rate}.

On a common probability space, the two norm-deviation bounds are
summable. Borel--Cantelli therefore gives
$c_0\sqrt m\le\|A_m\|_{\op}\le C\sqrt m$ eventually almost surely
for a deterministic $c_0>0$. The scaled logarithm of the norm then
converges almost surely to zero, and
Corollary~\ref{cor:almost-sure} supplies the almost-sure limit for
$L_{A_m}$. Finally, multiplication of $A$ by a positive scalar
multiplies both $U_A$ and $L_A$ by that scalar, leaving their ratio
unchanged.
\end{proof}

\end{document}